\documentclass{article}
\usepackage{graphicx}
\usepackage{amsmath}
\usepackage{epsfig}
\usepackage{latexsym}

\newtheorem{theorem}{Theorem}[section]

\newtheorem{lemma}[theorem]{Lemma}

\newenvironment{proof}[1][Proof]{\textbf{#1.} }{\ \rule{0.5em}{0.5em}}

\begin{document}

\title{Lower Bound for the Size of Maximal Nontraceable Graphs \thanks{%
This material is based upon work supported by the National Research
Foundation under Grant number 2053752.}}
\author{Marietjie Frick, Joy Singleton \\
University of South Africa,\\
P.O. Box 392, Unisa, 0003,\\
South Africa.\\
\textbf{e-mail:} singlje@unisa.ac.za frickm@unisa.ac.za}
\date{ }
\maketitle

\begin{abstract}
Let $g(n)$ denote the minimum  number of edges of a  maximal nontraceable graph of order $n$.
Dudek, Katona and Wojda (2003) showed that $ g(n) \geq \lceil \tfrac{3n-2}{2} \rceil -2$ for  $ n \ge 20$ 
and
$ g(n) \le  \lceil \tfrac{3n-2}{2} \rceil$ for $ n \geq 54$ as well as for  
$n \in I= \{22,23,30,31,38,39,40,41,42,43,46,47,48,49,50,51\}$. 
We show that $ g(n) =  \lceil \tfrac{3n-2}{2} \rceil$  for $ n \geq 54$ as well as for  
$n \in I \cup \{12,13\}$ and we determine $g(n)$ for $n \le 9$.


\end{abstract}

\section{\ Introduction}

We consider only simple, finite graphs $G$ and denote the vertex set, 
the edge set, the order and the size of $G$ by $V(G)$, $E(G)$, $v(G)$ and $e(G)$, respectively. 
The \textit{open neighbourhood} 
of a vertex $v$ in $G$ is the set 
$N_{G}(v)=\{x\in V(G):vx\in E(G)\}$. If $U$ is a nonempty subset of $V(G)$ 
then $\langle U \rangle$ denotes the subgraph of $G$ induced by $U$. 

 A graph $G$ is 
\textit{
hamiltonian }if it has a \textit{hamiltonian cycle }(a cycle containing all
the vertices of $G$), and 
\textit{traceable} if it has a \textit{hamiltonian
path }(a path containing all the vertices of $G$). 
A graph$\ G$ is \textit{
maximal nonhamiltonian}\emph{\ }(MNH)\ if $G$ is not hamiltonian, but $\ G+e$
is hamiltonian for each $e$ $\in E(\overline{G})$, where $\overline{G}$
denotes the complement of $G.$ 
A graph $G$ is \textit{maximal nontraceable}
(MNT) if $G$ is not traceable, but\ $G+e$ is traceable for each 
$e\in E(\overline{G})$.

\bigskip

In 1978 Bollob\'{a}s \cite{bollobas} posed the problem of finding the least
number of edges, $f(n)$, in a MNH graph of order $n$. Bondy \cite{bondy} had
already shown that a MNH graph with order $n\geq 7$ that contained $m$
vertices of degree $2$ had at least $(3n+m)/2$ edges, and hence $f(n)\geq
\left\lceil 3n/2\right\rceil $ for $n\geq 7$. Combined results of Clark, Entringer and Shapiro 
\cite{ce}, \cite{ces} and Lin, Jiang, Zhang and Yang \cite
{ljzy} show that $f(n)=\left\lceil 3n/2\right\rceil $ for $n\geq 19$ and
for $n=6,10,11,12,13,17$. The values of $f(n)$ for the remaining values of $n$ are also given in \cite
{ljzy}.

\bigskip

Let $g(n)$ denote the minimum number of edges in a MNT graph of order $n$. 
Dudek, Katona and Wojda \cite{dudek} proved that 
\[ g(n) \geq \lceil \tfrac{3n-2}{2} \rceil -2 \ \text{for } n \ge 20 \]
and showed, by construction, that
\[ g(n) \le  \lceil \tfrac{3n-2}{2} \rceil  \ \text{for } n \geq 54 \]  
\[ \text{as well as for }  
n \in I= \{22,23,30,31,38,39,40,41,42,43,46,47,48,49,50,51\}. \]
We prove, using a  method different from that in \cite{dudek}, that 
\[ g(n) \ge  \lceil \tfrac{3n-2}{2} \rceil  \ \text{for } n \geq 10. \] 
We also construct graphs of order $n=12,13$ with $ \lceil \frac{3n-2}{2} \rceil$ edges and thus show that 
\[ g(n) =  \lceil \tfrac{3n-2}{2} \rceil  \ \text{for } n \geq 54 \  \text {as well as for }  
n \in I \cup \{12,13\}. \]
We also determine $g(n)$ for $n \le 9$.

\section{Auxilliary Results}

In this section we present some results concerning MNT graphs,
which we shall  use, in the next section, to prove that a MNT graph of order $n \ge 10$ has 
at least $\frac{3n-2}{2}$ edges. The first one concerns the lower bound for the number of edges of MNH graphs. 
It is the combination of results proved in \cite{bondy} and \cite{ljzy}.

\begin{theorem}
\label{MNH}
(Bondy and Lin, Jiang, Zhang and Yang) If $G$ is a MNH graph of order $n$, then $e(G) 
\ge \frac{3n}{2}$ for $n \ge 6$.
\end{theorem}

The following lemma, which we proved in \cite{fs}, will be used frequently.

\begin{lemma}
\label{subgraph}Let $Q$ be a path in a MNT graph $G$. If $\langle V(Q) \rangle$ 
is not complete, then some internal vertex of $Q$ has a neighbour in $G-V(Q)$. 
\end{lemma}

\begin{proof}
Let $u$ and $v$ be two nonadjacent vertices of $Q$. Then $G+uv$ has a hamiltonian 
path $P$. Let $x$ and $y$ be the two endvertices of $Q$ and suppose no internal vertex 
of $Q$ has a neighbour in $G-V(Q)$. Then $P$ has a subpath $R$ in 
$\langle V(Q) \rangle + uv$ and $R$ has either one or both endvertices in $\{x,y\}$. If 
$R$ has only one endvertex in $\{x,y\}$, then $P$ has an endvertex in $Q$. In either 
case the path obtained from $P$ by replacing $R$ with $Q$ is a hamiltonian path of $G$. 
\hfill \end{proof}

\bigskip

The following lemma is easy to prove.

\begin{lemma}
\label{cutset}
Suppose $T$ is a cutset of a connected graph $G$ and  $A_1,...,A_k$ are components of $G-T$.

(a) If $k \ge |T|+2$, then $G$ is nontraceable.

(b) If $G$ is MNT then $k \le |T|+2$.

(c) If $G$ is MNT and $k=|T|+2$, then $\langle T \cup A_i \rangle$ is complete 
for $i=1,2,...,k$.

\end{lemma}

\begin{proof}
(a) and (b) are obvious. If (c) is not true, then there is an $i$ such that 
$\langle T \cup A_i \rangle$ has two nonadjacent vertices $x$ and $y$. But then $T$ is 
a cutset of the graph $G+xy$ and $(G+xy)-T$ has $|T|+2$ components and hence $G+xy$ is nontraceable, by (a).

\hfill \end{proof}

The proof of the following lemma is similar to the previous one.

\begin{lemma}
\label{block} Suppose $B$ is a block of a connected graph $G$.

(a) If $B$ has more than two cut-vertices, then $G$ is nontraceable.

(b) If $G$ is MNT, then $B$ has at most three cut-vertices.

(c) If $G$ is MNT and $B$ has exactly three cut-vertices, then $B$ consists of exactly 
four blocks, each of which is complete.

\end{lemma}

In \cite{fs} we proved some results concerning the degrees of the neighbours 
of the vertices of degree 2 in a 2-connected MNT graph, which enabled us to show that the 
average degree of the vertices in a 2-connected MNT graph is at least 3. We now restate those 
results in a form that is applicable also to MNT graphs which are not 2-connected. (Note that 
in a 2-connected graph no two vertices of degree 2 are adjacent to one another.)

 \bigskip 

\begin{lemma}
\label{mntdeg2}If $G$ is a connected MNT graph and $v\in V(G)$ with $d\left( v\right)
=2$, then the neighbours of $v$ are adjacent. Also, one of the neighbours has degree at least $4$ 
and the other neighbour has degree $2$ or at least $4$.
\end{lemma}

\begin{proof}
Let $N_G(v)=\{x_1,x_2\}$ and let $Q$ be the path $x_1vx_2$. Since $N_G(v)\subseteq Q$, 
it follows from Lemma \ref{subgraph} that $\langle V(Q) \rangle$ is a complete graph;
 hence $x_1$ and $x_2$ are adjacent.

Since $G$ is connected and nontraceable, at least one of $x_1$ and $x_2$ has degree bigger that 2. 
Suppose $d(x_{1})>2$ and let $z \in N(x_1)- \{v,x_2\}$. If $Q$ is the path $zx_1vx_2$ then, since $d(v)=2$, 
the graph $\langle V(Q) \rangle$ is not complete and hence it follows from 
Lemma \ref{subgraph} that $d(x_1)\geq 4$. Similarily if $d(x_2) >2$, then $d(x_2)\geq 4$ .
\hfill \end{proof}

\bigskip

\begin{lemma}
\label{mnt1deg2}Suppose $G$ is a connected MNT graph with distinct nonadjacent vertices $v_{1}$ and 
$v_{2}$ such that $d(v_1)=d(v_2) =2$.

(a) If $v_{1}$ and $v_{2}$ have exactly one common neighbour $x$, then $d(x)\geq 5.$

(b) If $v_{1}$ and $v_{2}$ have the same two neighbours $x_{1}$ and $x_{2}$, then $N_{G}(x_{1})-\{x_{2}
\}=N_{G}(x_{2})-\{x_{1}\}$ and  $d(x_{1})=d(x_{2})\geq 5.$
\end{lemma}

\begin{proof}
(a) Let $N(v_i)=\{x,y_i\}$; $i=1,2$. It follows from Lemma \ref{mntdeg2} that 
$x$ is adjacent to $y_i$; $i=1,2$. Let $Q$ be the path $y_1v_1xv_2y_2$. Since 
$\langle V(Q) \rangle$ is not complete, it follows from Lemma \ref{subgraph} that $x$ has 
a neighbour in $G-V(Q)$. Hence $d(x) \geq 5$.

(b) From Lemma \ref{mntdeg2} it follows that $x_{1}$ and $x_{2}$ are adjacent.
Let $Q$ be the path $x_2v_1x_1v_2$. $\langle V(Q) \rangle$ is not complete since $v_1$ and $v_2$ 
are nonadjacent. Thus it follows from Lemma \ref{subgraph} that $x_1$ has a neighbour in $G-V(Q)$. Now 
suppose $p \in N_{G-V(Q)}(x_1)$ and  $p \notin N_{G}(x_2)$. 
Then a hamiltonian path $P$ in $G+px_{2}$ contains a
subpath of either of the forms given in the first column of Table 1. Note
that $i,j\in \{1,2\}$; $i\neq j$ and that $L$ represents a subpath of $P$ in $
G-\{x_{1},x_{2},v_{1},v_{2},p\}$. If each of the subpaths is replaced by the
corresponding subpath in the second column of the table we obtain a
hamiltonian path $P^{\prime }$ in $G$, which leads to a contradiction.

\begin{equation*}
\begin{tabular}{|l|l|}
\hline
Subpath of $P$ & Replace with \\ \hline
$v_{i}x_{1}v_{j}x_{2}p$ & $v_{i}x_{2}v_{j}x_{1}p$ \\ \hline
$v_{i}x_{1}Lpx_{2}v_{j}$ & $v_{i}x_{2}v_{j}x_{1}Lp$ \\ \hline
\end{tabular}
\end{equation*}
\begin{equation*}
\text{Table 1}
\end{equation*}
Hence $p \in N_{G}(x_2)$. Thus $N_{G}(x_{1})-\{x_{2}\}\subseteq N_{G}(x_{2})-\{x_{1}\}$. 
Similarly $N_{G}(x_{2})-\{x_{1}\}\subseteq N_{G}(x_{1})-\{x_{2}\}$. Thus $
N_{G}(x_{1})-\{x_{2}\}=N_{G}(x_{2})-\{x_{1}\}$ and hence $d(x_{1})=d(x_{2})$. Now let 
$Q$ be the path $px_1v_1x_2v_2$. Since $\langle V(Q) \rangle$ is not complete, 
it follows from Lemma \ref{subgraph} that $x_1$ or $x_2$ has a neighbour in $G-V(Q)$. 
Hence $d(x_1)=d(x_2) \geq 5$.
\hfill \end{proof}

\begin{lemma}
\label{mnt3deg2}Suppose $G$ is a connected MNT graph of order $n\geq 6$ and that $v_{1},v_{2}$ and $
v_{3}$ are vertices of degree $2$ in $G$ having the same neighbours, $x_{1}$
and $x_{2}$. Then $G-\{v_{1},v_{2},v_{3}\}$ is complete and hence $e(G)=
\frac{1}{2}(n^{2}-7n+24)$.
\end{lemma}

\begin{proof}
The set $\{x_1,x_2\}$ is a cutset of $G$. Thus according to Lemma \ref{cutset} 
$G-\{v_{1},v_{2},v_{3}\}=K_{n-3}.$ Hence $e(G)=\frac{1}{2}(n-3)(n-4)+6.$
\hfill \end{proof}

\bigskip

By combining the previous three results we obtain

\begin{theorem}
\label{size1}
Suppose $G$ is a connected  MNT graph without vertices of degree $1$ or adjacent vertices of degree $2$. 
If $G$ has order $n\geq 7$ and $m$ vertices of degree $2$, then $e(G)\geq \frac{1}{2}(3n+m)$.
\end{theorem}

\begin{proof}
If $G$ has three vertices of degree 2 having the same two neighbours then, by 
Lemma \ref{mnt3deg2}, $m=3$ and
\begin{equation*}
e(G)=\tfrac{1}{2}(n^{2}-7n+24)\geq \tfrac{1}{2}(3n+m)\ \text{when}\ n \geq 7.
\end{equation*}

We now assume that $G$ does not have three vertices of degree 2 that have
the same two neighbours. Let $v_{1},...,v_{m}$ be the vertices of degree 2
in $G$ and let $H=G-\{v_{1},...,v_{m}\}.$ Then by Lemmas \ref{mntdeg2} and  \ref
{mnt1deg2} the minimum degree, $\delta (H)$ of $H$ is at
least 3. Hence 
\begin{equation*}
e(G)=e(H)+2m\geq \tfrac{3}{2}(n-m)+2m=\tfrac{1}{2}(3n+m).
\end{equation*}
\vspace{-2mm} \hfill \end{proof}

\bigskip

\section{The minimum size of a MNT graph}

Our aim is to determine the exact value of $g(n)$.
By consulting the Atlas of Graphs \cite{readwilson}, one can see, by inspection, that $g(2)=0$, 
$g(3)=1$, $g(4)=2$, $g(5)=4$, $g(6)=6$ and $g(7)=8$ (see Fig.\ 3).

We now give a lower bound for $g(n)$ for $n \ge 8$.

\begin{theorem}
\label{minsize}
If $G$ is a MNT graph of order $n$, then
\begin{eqnarray*}
e(G) &\geq &\left\{ 
\begin{tabular}{ll}
$10$ &if $n=8$   \\
$12$ &if $n=9 $  \\ 
$\frac{3n-2}{2}$ & if $ n \ge 10.$   
\end{tabular} 
\right. \\
\end{eqnarray*}

\end{theorem}

\begin {proof}

If $G$ is not connected, then $G=K_k \cup K_{n-k}$, for some positive integer $k<n$ and then, clearly, 
$e(G) > \frac{3n-2}{2}$ for $n \ge 8$. Thus we assume that $G$ is connected.

We need to prove that the sum of the degrees of the vertices of $G$ is at least $3n-2$. In view of 
Theorem \ref{size1}, we let
\[M= \{v \in V(G) \  |\  d(v)=2\  \text{and no neighbour of } v\  \text{has degree } 2 \}. \]
The remaining vertices of degree 2 can be dealt with simultaneously with the vertices of degree 1. We let
\[S=\{v \in V(G)-M\  |\  d(v)=2\  \text{or } d(v)=1 \}.\]

If $S= \emptyset$, then it follows from Theorem \ref{size1} that $e(G) \ge \frac{1}{2}(3n+m)$. Thus 
we assume that $S \neq \emptyset$.

We observe that, if $H$ is a component of the graph of $\langle S \rangle$, then either  
$H \cong K_1$ or $ H \cong K_2$ and $N_G(H)-V(H)$ consists of a single vertex, which is a cut-vertex of $G$.

An example of such a graph $G$ is depicted in the figure below.
\[
\includegraphics{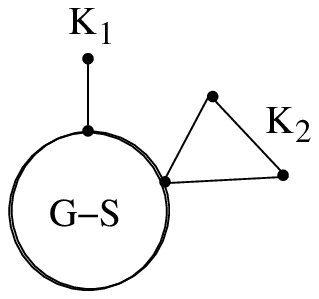}
\]
\[
\text{Fig.\ 1}
\]

Let $s=|S|$. By Lemma \ref{block} the graph $ \langle S \rangle $ has at most three components. We thus have three cases:

\newpage

\textbf{CASE 1. } $\langle S \rangle$ has exactly three components, say $H_1,H_2,H_3$:

In this case the neighbourhoods of $H_1,H_2,H_3$ are pairwise disjoint; hence $G$ has three cut-vertices. Hence it follows 
from Lemma \ref{block} that $G-S$ is a complete graph of order at least 3. Futhermore, for every possible value of $s$, 
the number of edges in $G$ incident with the vertices in $S$ is $2s-3$. Thus
\[e(G)=  \binom{n-s}{2} +2s -3 \ \text{for } s=3,4,5 \ \text{or } 6;\ s \le n-3. \]
An easy calculation shows that, for each possible value of $s$,
\begin{eqnarray*}
e(G) &\geq &\left\{ 
\begin{tabular}{ll}
$10$ &if $n=8$   \\
$12$ &if $n=9 $  \\ 
$\frac{3n-2}{2}$ & if $ n \ge 10.$   
\end{tabular} 
\right. \\
\end{eqnarray*}

This case is a Zelinka Type II construction, cf.\ \cite{zelinka}. 
The graphs of smallest size of order 8 and 9 given by this construction are depicted in Fig.\ 3.

\bigskip

\textbf{CASE 2. }  $\langle S \rangle$ has exactly two components, say $H_1,H_2$:

In this case the number of edges in $G$ incident with the vertices in $S$ is $2s-2$.

\bigskip

\textbf{Subcase 2.1. } $N_{G}(H_1)=N_{G}(H_2)$:\\
Then it follows from Lemma \ref{cutset} that $G-S$ is a complete graph.  Hence
\[ e(G)= \binom{n-s}{2} +2s -2 \ \text{for } s=2,3 \ \text{or } 4.  \]
Thus
\begin{eqnarray*}
e(G) &\geq &\left\{ 
\begin{tabular}{ll}
$12$ &if $n=8$   \\
$16$ &if $n=9 $  \\ 
$\frac{3n-2}{2}$ & if $ n \ge 10$   
\end{tabular} 
\right. \\
\end{eqnarray*}

This case is a Zelinka Type I construction, cf.\ \cite{zelinka}.

\bigskip

\textbf{Subcase 2.2. } $N_{G}(H_1) \neq N_{G}(H_2)$:\\
Let  $N_{G}(H_i)=y_i$, $i=1,2$ and $y_1 \neq y_2$. \\
If $y_1y_2 \notin E(G)$ then $G+y_1y_2$ has a hamiltonian path $P$. 
But then $P$ has one endvertex in $H_1$ and the other in $H_2$ and contains the edge  
$y_1y_2$; hence $V(G-S) =\{y_1,y_2 \}$. But then $G$ is disconnected. 
This contradiction shows that  $y_1y_2 \in E(G)$. 

Now $G-S$ is not complete, otherwise $G$ would be traceable. Since $G+vw$, where 
$v$ and $w$ are nonajacent vertices in $V(G-S)$,  contains a hamiltonian path with one 
endvertex in $H_1$ and the other in $H_2$ and $y_1y_2 \in E(G)$, it follows that $(G-S)+vw$ 
has a hamiltonian cycle. Hence $G-S$ is either hamiltonian or MNH. We consider  these two cases separately:

\bigskip 

\textbf{Subcase 2.2.1. } $G-S$ is hamiltonian: \\
Then no hamiltonian cycle in $G-S$ contains $y_1y_2$, otherwise $G$ would be traceable. 
Thus $d_{G-S}(y_i) \ge 3$ for $i=1,2$.

It also follows from Lemma \ref{cutset} that no vertex $v \in M$ can be adjacent to both $y_1$ and $y_2$ 
since the graph $\langle V(H_i) \cup T \rangle $, where $T=\{y_1,y_2\}$ is not complete, for $i=1,2$.
If $v \in M$ is adjacent to to one of the $y_i$'s for $i=1,2$, say $y_1$, then, since the neighbours of 
$v$ are adjacent, it follows that $d_{G-M- S}(y_1) \ge 3$.

It follows from our definition of $M$ and $S$ that $N_G(M) \cap S = \emptyset$. Since $G-M$ is not a 
complete graph, it follows from Lemma \ref{mnt3deg2} that $M$ does not have three vertices that have the 
same neighbourhood in $G$. Hence, by Lemmas \ref{mntdeg2} and \ref{mnt1deg2}, the minimum degee of 
the graph $G- M- S$ is at least 3.

 Now, for $n \ge 8$
\begin{eqnarray*}
e(G) &=&  e(G-M-S)+2m+2s-2 \\
&\geq &\frac{1}{2}\left(3\left(n-m-s\right)\right)+2m+2s-2   \\
&=&\frac{1}{2}\left(3n+m+s-4\right) \\
& \ge &  \frac{3n-2}{2}, \ \text{since } s \ge 2. 
\end{eqnarray*}

\textbf{Subcase 2.2.2. } $G-S$ is nonhamiltonian: \\
Then  $G-S$ is MNH (as shown above); hence it follows from Theorem \ref{MNH}, that \\
$e(G-S) \ge \frac{3}{2}(n-s)$ for $n-s \ge 6$. 

Thus, for  $n-s \ge 6$ and $n \ge 8$
\begin{eqnarray*}
e(G)& =& e(G-S)+2s-2 \\
&\ge & \frac{1}{2}(3(n-s)) +2s-2 \\
&=& \frac{1}{2}(3n+s-4) \\
&\ge & \frac{3n-2}{2}, \ \text{since } s \ge 2. 
\end{eqnarray*}

From \cite{ljzy} we have 
\begin{eqnarray*}
e(G-S) &\geq &\left\{ 
\begin{tabular}{ll}
$6$ &for $n-s=5$   \\
$4$ &for $n-s=4 $.   
\end{tabular} 
\right. \\
\end{eqnarray*}
Thus
\begin{eqnarray*}
e(G) &\geq &\left\{ 
\begin{tabular}{ll}
$12$ &for $n=9$ and $n-s=5$   \\
$10$ &for $n=8$ and  $n-s=5 $ or $n-s=4$. 
\end{tabular} 
\right. \\
\end{eqnarray*} 

The smallest MNH graphs $F_4$ and $F_5$ of order 4 and 5 respectively, are depicted in Fig.\ 2; cf.\ \cite {ljzy}. 
The graphs $G_8$ and $G_9$ (see Fig.\ 3)  are obtained, respectively, by using $F_4$ with $s=4$ or $F_5$ with $s=3$, 
and $F_5$ with $s=4$.
\[
\includegraphics{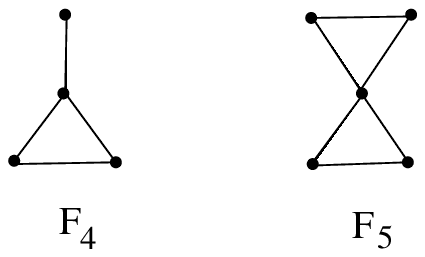}
\]
\[
\text{Fig.\ 2}
\]
\textbf{CASE 3. } $\langle S \rangle$ has exactly one component, say $H$: \\
Since
\[ \sum_{v \in S}d_{G}(v) = 3s-2, \ \text{for }s=1,2 \]
it follows that
\begin{eqnarray*}
e(G) &=&  e(G-M)+2m \\
&=&\frac{1}{2} \left( \sum_{v \in V(G-M)-S}d_{G-M}(v) + \sum_{v \in S}d_{G-M}(v) \right) + 2m \\
&\geq & \frac{1}{2}\left(3\left(n-m-s\right)+3s-2\right)+2m \\
&=&\frac{1}{2}\left(3n+m-2\right) \\
&\ge & \frac{3n-2}{2}.
\end{eqnarray*}

\vspace{-4mm}\hfill \end{proof}

From the previous theorem we have $g(8)=10$, $g(9)=12$ and $g(n) \ge \lceil \frac{3n-2}{2} \rceil$ for $n \ge 10$. 
The MNT graphs $G_n$ of order $n$ with $g(n)$ edges, for $n \le 9$ are given in Fig.\ 3.
\[
\includegraphics{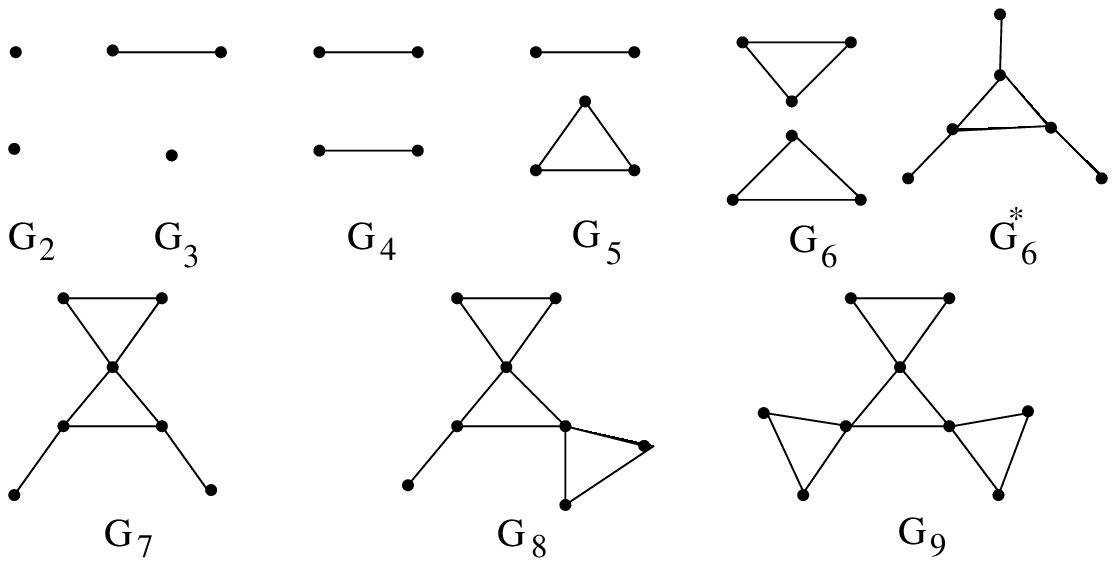} 
\]
\[
\text{Fig.\ 3}
\]

In \cite{dudek} Dudek, Katona and Wojda constructed, for every 
$n \geq 54$ as well as for every $n \in I= \{22,23,30,31,38,39,40,41,42,43,46,47,48,49,50,51\}$,
 a MNT graph of size $\lceil \frac{3n-2}{2} \rceil$ in the following way:
Consider a cubic MNH graph $G$  with the property that \\
(1) there is an edge $y_1y_2$ of
 $G$, such that $N(y_1) \cap N(y_2) =  \emptyset $, and  \\
(2) $G+e$ has a hamiltonian 
cycle containing $y_1y_2$ for every $e \in E(\overline G) $. 

Now take two graphs $H_1$ and $H_2$, with $H_1 \cong K_1$ and $H_2 \cong K_1$ or $H_2 \cong K_2$ and join each vertex of 
$H_i$ to $y_i;\ i=1,2$. The new graph is a MNT graph of order $v(G)+2$ and size $e(G)+2$ or of order $v(G)+3$ and 
size $e(G)+4$. 

It follows from results in \cite{ce} and \cite{ces} that for every even $n \ge 52$ as well as for 
$n \in \{20,28,36,38,40,44,46,48\}$ there exists a cubic MNH  graph of order $n$ that satisfies (1) and (2). 
Thus this construction provides MNT graphs of order $n$ and size $\lceil \frac{3n-2}{2} \rceil$ for every $n \ge 54$ 
as well as for every  $n \in I $.

We determined, by using the Graph Manipulation Package developed by Siqinfu and Sheng Bau*, 
that the Petersen graph  also satisfies the above 
property. Hence, according to the above construction, there are also MNT graphs of order $n$ and size 
$\lceil  \frac{3n-2}{2} \rceil $ for $n=12,13$.

Thus $g(n) = \lceil  \frac{3n-2}{2} \rceil$ for $n \geq 54$ as well as for every $n \in I \cup \{12,13\}$.

It remains an open problem to find $g(n)$ for $n=10,11$ and those values of $n$ between 13 and 54 which are not in $I$.
 
\bigskip

*\textbf{Acknowledgement} We wish to thank Sheng Bau for allowing us the use of the 
programme, Graph Manipulation Package  Version 1.0 (1996), Siqinfu and Sheng Bau,
Inner Mongolia Institute of Finance and Economics, Huhhot, CN-010051, People's 
Republic of China.

\end{document}